\documentclass[12pt,leqno]{amsart}
\newtheorem{dummy}{}[section]

\newtheorem{theorem}[dummy]{Theorem}

\newtheorem{remark}[dummy]{Remark}

\begin{document}
\bibliographystyle{plain}
\title{Non-factorisation of Arf-Kervaire classes through ${\mathbb RP}^{\infty} \wedge  {\mathbb RP}^{\infty}$}
\author{Victor P. Snaith}
\date{20 January 2010}
\maketitle
 
 \begin{abstract}
As an application of the upper triangular technology method of \cite{Arfbk08} it is shown that there do not exist stable homotopy classes of  $ {\mathbb RP}^{\infty} \wedge  {\mathbb RP}^{\infty}$ in dimension $2^{s+1}-2$ with $s \geq 2$ whose composition with the Hopf map to $ {\mathbb RP}^{\infty}$ followed by the Kahn-Priddy map gives an element in the stable homotopy of spheres of Arf-Kervaire invariant one.
 \end{abstract}
 
 \section{Introduction}
 \begin{dummy}
\label{1.1}
\begin{em}

For $n > 0 $ let $\pi_{n}(\Sigma^{\infty}S^{0})$ denote the $n$-th stable homotopy group of $S^{0}$, the $0$-dimensional sphere.
Via the Pontrjagin-Thom construction an element of this group corresponds to a framed bordism class of an $n$-dimensional framed manifold. The Arf-Kervaire invariant problem concerns whether or not there exists such a framed manifold possessing a Kervaire surgery invariant which is non-zero (modulo $2$). In \cite{Brow69} it is shown that this can happen only when $n = 2^{s+1}-2$ for some $s \geq 1$. Resolving this existence problem is an important unsolved problems in homotopy theory (see \cite{Arfbk08} for a historical account of the problem together with new proofs of all that was known up to 2008). Recently important progress has made (\cite{HHR09}; see also \cite{Akh08}, \cite{Akh08a}) which shows that $n = 126$ is the only remaining possibility for existence (more details may be found in the survey article \cite{Arfsurv10}.

In view of the renewed interest in the Arf-Kervaire invariant problem it may be of interest to describe a related non-existence result. An equivalence formulation (see \cite{Arfbk08} \S1.8) is that there exists a stable homotopy class
$\Theta : \Sigma^{\infty} S^{2^{s+1}-2} \longrightarrow  \Sigma^{\infty} {\mathbb RP}^{\infty} $ with mapping cone ${\rm Cone}(\Theta)$ such that the Steenrod operation
\[ Sq^{2^{s}} : H^{2^{s} - 1}( {\rm Cone}(\Theta) ; {\mathbb Z}/2)  \cong {\mathbb Z}/2 \longrightarrow
H^{2^{s+1} - 1}( {\rm Cone}(\Theta) ; {\mathbb Z}/2)  \]
is non-trivial. Using the upper triangular technology (UTT) of \cite{Arfbk08} we shall prove the following result:
\end{em}
\end{dummy}
\begin{theorem}{$_{}$}
\label{1.2}
\begin{em}

Let $H: \Sigma^{\infty} {\mathbb RP}^{\infty} \wedge  {\mathbb RP}^{\infty}  \longrightarrow  
\Sigma^{\infty} {\mathbb RP}^{\infty} $ denote the map obtained by applying the Hopf construction to the multiplication on ${\mathbb RP}^{\infty}$. Then, if $s \geq 2$, there does not exist a stable homotopy class
\[  \tilde{\Theta} : \Sigma^{\infty} S^{2^{s+1}-2} \longrightarrow  \Sigma^{\infty} {\mathbb RP}^{\infty} \wedge  {\mathbb RP}^{\infty} \] 
such that the composition $\Theta = H \cdot \tilde{\Theta}$ is detected by a non-trivial $Sq^{2^{s}}$ as in \S\ref{1.1}.
\end{em}
\end{theorem}

In  \S\ref{2.2} this result will be derived as a simple consequence of the UTT relations (\cite{Arfbk08} Chapter Eight). The basics of the UTT method are sketched in \S\ref{2.1}. Doubtless there are other ways to prove Theorem \ref{1.2} (for example, from the results of \cite{ST82}; see also \cite{Arfbk08} Chapter  Two) but it provides an elegant application of UTT.

 \section{Upper triangular technology (UTT)}
\begin{dummy}
\label{2.1}
\begin{em}

Let $F_{2n}( \Omega^{2} S^{3} )$ denote the $2n$-th filtration of the combinatorial model
for $\Omega^{2} S^{3} \simeq W \times S^{1}$. Let $F_{2n}( W )$
denote the induced filtration on $W$ and let $B(n)$ be the Thom spectrum of the 
canonical bundle induced by $f_{n} :  \Omega^{2} S^{3} \longrightarrow BO$, where $B(0) = S^{0}$ by convention. From \cite{Mah81} one has a $2$-local, left $bu$-module homotopy equivalence of the form\footnote{In \cite{Arfbk08} and related papers I consistently forgot what I had written in my 1998 McMaster University notes ``On $bu_{*}(BD_{8})$''. Namely, in the description of Mahowald's result I stated that $\Sigma^{4n}B(n) $ was equal to the decomposition factor $F_{4n}/F_{4n-1}$ in the Snaith splitting of $ \Omega^{2} S^{3} $. Although this is rather embarrassing, I got the homology correct so that the results remain correct upon replacing $F_{4n}/F_{4n-1}$ by $\Sigma^{4n}B(n) $ throughout! I have seen errors like this in the World Snooker Championship where the no.$1$ player misses an easy pot by concentrating on positioning the cue-ball. In mathematics such errors are inexcusable whereas in snooker they only cost one the World Championship.}
\[   \bigvee_{n \geq 0 }  bu  \wedge  \Sigma^{4n}B(n)  
\stackrel{ \simeq }{\longrightarrow}  bu  \wedge bo .\]
Therefore, if $ \Theta$ is as in \S\ref{1.1}, then 
\[(bu  \wedge bo)_{*}({\rm Cone}(\Theta))  \cong  \bigoplus_{n \geq 0} \  (bu_{*}({\rm Cone}(\Theta) \wedge \Sigma^{4n}B(n) ) .  \]

Let $\alpha(k)$ denote the number of $1$'s in the dyadic expansion of the positive integer $k$. For $1 \leq k \leq 2^{s-1}-1$ and $2^{s} \geq 4k - \alpha(k) + 1$ there are isomorphisms of the form (\cite{Arfbk08} Chapter Eight \S4)
\[   \begin{array}{l}
    bu_{2^{s+1}-1}( C(\Theta)  \wedge \Sigma^{4k}B(k) )   \cong bu_{2^{s+1}-1}( {\mathbb RP}^{\infty}  \wedge \Sigma^{4k}B(k) )      \cong  V_{k} \oplus 
 {\mathbb Z}/2^{2^{s}-4k + \alpha(k) }
\end{array}   \]
where $V_{k}$ is a finite-dimensional ${\mathbb F}_{2}$-vector space consisting of elements which are detected in mod $2$ cohomology (i.e. in filtration zero, represented on the $s=0$ line) in the mod $2$ Adams spectral sequence. The map $1 \wedge \psi^{3} \wedge 1$ on $bu \wedge bo \wedge C(\Theta)$ acts on the direct sum decomposition like the upper triangular matrix
\[   \left(\begin{smallmatrix}1 & 1 & 0 & 0 & 0 & \ldots\\0 & 9 & 1 & 0 & 0 & \ldots\\0 & 0 & 9^{2} & 1 & 0 Ê& \ldots\\0 & 0 & 0 & 9^{3} & 1 Ê & \ldots \\\vdots & Ê\vdots & Ê\vdots & Ê\vdots & Ê\vdots & Ê\vdots \end{smallmatrix}\right) . \]
In other words $(1 \wedge \psi^{3} \wedge 1)_{*}$ sends the $k$-th summand to itself by multiplication by $9^{k-1}$ and sends the $(k-1)$-th summand to the $(k-2)$-th by a map  
\[    \begin{array}{l}
 (\iota_{k,k-1})_{*}   :   V_{k} \oplus 
 {\mathbb Z}/2^{2^{s}-4k + \alpha(k) }   \longrightarrow  V_{k-1} \oplus 
 {\mathbb Z}/2^{2^{s}-4k +4+ \alpha(k-1) }   
 \end{array}   \]
 for $2 \leq k \leq 2^{s-1}-1$ and $2^{s} \geq 4k - \alpha(k) + 1$. The right-hand component of this map is injective on the summand $  {\mathbb Z}/2^{2^{s}-4k + \alpha(k) } $ and annihilates $V_{k}$. 
 
 It is shown in \cite{Knapp97} (also proved by UTT in (\cite{Arfbk08} Chapter Eight when $s \geq 2$) that $\Theta$ corresponds to a stable homotopy class of Arf-kervaire invariant one if and only if it is detected by the Adams operation $\psi^{3}$ on $\iota \in bu_{2^{s+1}-1}({\rm Cone}(\Theta))$, an element of infinite order.
 
From these properties and the formula for $\psi^{3}(\iota)$ one easily obtains a series of equations (\cite{Arfbk08} \S8.4.3) for the components of $(\eta \wedge 1 \wedge 1)_{*}(\iota)$ where $\eta : S^{0} \longrightarrow bu$ is the unit of $bu$-spectrum. Here we have used the isomorphism $bu_{2^{s+1}-1}(C(\Theta)) \cong bo_{2^{s+1}-1}(C(\Theta)) $ since, strictly speaking,  the latter group is the domain of 
$(\eta \wedge 1 \wedge 1)_{*}$.  It is shown in (\cite{Arfbk08} Theorem 8.4.7) that this series of equations implies that the 
$bu_{2^{s+1}-1}({\rm Cone}(\Theta) \wedge \Sigma^{2^{s}}B(2^{s-2}) ) $-component of 
$(\eta \wedge 1 \wedge 1)_{*}(\iota)$ is non-trivial and gives some information on the identity of this non-trivial element.

It is this information which we shall now use to prove Theorem \ref{1.2}.
\end{em}
\end{dummy}
\begin{dummy}{Proof of Theorem \ref{1.2}}
\label{2.2}
\begin{em}

Suppose, for a contradiction,  that $\Theta$ and $\tilde{\Theta}$ exist. We must assume that $s \geq 2$ because the UTT results of (\cite{Arfbk08} Theorem 8.4.7) are only claimed for this range.

The mod $2$ cohomology of $ \Sigma^{2^{s}}B(2^{s-2}) $ is given by the ${\mathbb F}_{2}$-vector space with basis $\{z_{2^{s} + 2j}, 0 \leq j \leq  2^{s-1}-2;
z_{2^{s}+3 + 2k}, 0 \leq k \leq  2^{s-1}-2 \}$ on which the left action by $Sq^{1}$ and $Sq^{0,1} = Sq^{1}Sq^{2} + Sq^{2}Sq^{1}$
are given by $Sq^{1}(z_{2^{s} + 2j}) = z_{2^{s} + 3 + 2(j-1)}$ for $1 \leq j \leq 2^{s-1}-1$ and $Sq^{0,1}(z_{2^{s} + 2j}) = z_{2^{s} + 3 + 2j}$ for $0 \leq j \leq 2^{s-1}-2$ and $Sq^{1}, Sq^{0,1}$ are zero otherwise. This cohomology module is the ${\mathbb F}_{2}$-dual of the ``lightning flash'' module depicted in (\cite{Ad74} p.341). 

Now consider the two $2$-local Adams spectral sequences
\[  \begin{array}{l}
E_{2}^{s,t} =   {\rm Ext}_{B}^{s, t}( H^{*}( C(\Theta)   ; {\mathbb Z}/2) \otimes  H^{*}(  \Sigma^{2^{s}}B(2^{s-2}; {\mathbb Z}/2) , {\mathbb Z}/2)  \\
\\
\hspace{90pt}  \Longrightarrow  bu_{t-s}(C(\Theta)  \wedge \Sigma^{2^{s}}B(2^{s-2}) ),
\end{array}  \]
which collapses and
\[  \begin{array}{l}
\tilde{E}_{2}^{s,t} =   {\rm Ext}_{B}^{s, t}( H^{*}( C(\tilde{\Theta})   ; {\mathbb Z}/2) \otimes  H^{*}(  \Sigma^{2^{s}}B(2^{s-2}; {\mathbb Z}/2) , {\mathbb Z}/2)  \\
\\
\hspace{90pt}  \Longrightarrow  bu_{t-s}(C(\tilde{\Theta})  \wedge \Sigma^{2^{s}}B(2^{s-2}) ),
\end{array}  \]
where $B$ is the exterior subalgebra of the mod $2$ Steenrod algebra generated by $Sq^{1}$ and $ Sq^{0,1}$.

To fit in with the notation of  (\cite{Arfbk08} Theorem 8.4.7) set $s = q+2$ in Theorem \ref{1.2}. As mentioned in \S\ref{2.1}, it is shown in (\cite{Arfbk08} Theorem 8.4.7) that the component of $(\eta \wedge 1 \wedge 1)_{*}(\iota)$ lying in 
\[   \begin{array}{l}
bu_{2^{q+3}-1}(C(\Theta)  \wedge \Sigma^{2^{s}}B(2^{s-2}) )  \\
\\
 \cong bu_{2^{q+3}-1}({\mathbb RP}^{\infty}  \wedge \Sigma^{2^{s}}B(2^{s-2}) ) \\
 \\
  \cong   {\rm Ext}_{B}^{0, 2^{q+3}-1}( H^{*}( {\mathbb RP}^{\infty}   ; {\mathbb Z}/2) \otimes  H^{*}(  \Sigma^{2^{s}}B(2^{s-2}; {\mathbb Z}/2) , {\mathbb Z}/2) \\
  \\
  \subseteq  {\rm Hom}( \oplus_{u+v =  2^{q+3}-1}  \  H^{u}( {\mathbb RP}^{\infty}   ; {\mathbb Z}/2) \otimes  H^{v}(  \Sigma^{2^{s}}B(2^{s-2}; {\mathbb Z}/2) , {\mathbb Z}/2)
\end{array}  \]
corresponds to a homomorphism $f$ such that $f(x^{2^{q+2}-1} \otimes z_{2^{q+2}}) $ is non-trivial.

The factorisation $\Theta = H \cdot \tilde{\Theta}$ implies that there exists $h \in \tilde{E}_{\infty}^{0, 2^{q+3}-1} \subseteq 
\tilde{E}_{2}^{0, 2^{q+3}-1} $ such that $H_{*}(h)=f$. On the other hand
\[ \tilde{E}_{2}^{0, 2^{q+3}-1}  \cong   {\rm Ext}_{B}^{0, 2^{q+3}-1}( H^{*}( {\mathbb RP}^{\infty} \wedge {\mathbb RP}^{\infty}    ; {\mathbb Z}/2) \otimes  H^{*}(  \Sigma^{2^{s}}B(2^{s-2}; {\mathbb Z}/2) , {\mathbb Z}/2) .  \]
Therefore the homomorphism
The homomorphism
\[  \begin{array}{c}
 {\rm Ext}_{B}^{0, 2^{q+3}-1}( H^{*}( {\mathbb RP}^{\infty}  \wedge {\mathbb RP}^{\infty} ; {\mathbb Z}/2) \otimes  H^{*}( F_{2^{q+2}}/ F_{2^{q+2}-1} ; {\mathbb Z}/2) , {\mathbb Z}/2)    \\
 \\
 (H \wedge 1 )_{*}   \    \downarrow   \\
 \\
 {\rm Ext}_{B}^{0, 2^{q+3}-1}( H^{*}( {\mathbb RP}^{\infty}   ; {\mathbb Z}/2) \otimes  H^{*}( F_{2^{q+2}}/ F_{2^{q+2}-1} ; {\mathbb Z}/2) , {\mathbb Z}/2)
\end{array}    \]
satisfies $ (H \wedge 1 )_{*} (h)(x^{2^{q+2}-1} \otimes z_{2^{q+2}} )   = f(x^{2^{q+2}-1} \otimes z_{2^{q+2}})  \not\equiv 0$.
However
\[ \begin{array}{l}
 (H \wedge 1 )_{*} (h)(x^{2^{q+2}-1} \otimes z_{2^{q+2}} )   \\
 \\
 =  h(  \sum_{a=1}^{2^{q+2}-2}  \    x^{a} \otimes  x^{2^{q+2}-  a -1}  \otimes z_{2^{q+2}}  )  .
\end{array} \]
On the other hand
\[   \begin{array}{l}
Sq^{1}( x^{\alpha} \otimes  x^{2^{q+2}-2 - \alpha} \otimes z_{2^{q+2}}  ) \\
\\
=  \alpha ( x^{\alpha} \otimes  x^{2^{q+2}-1 - \alpha} \otimes z_{2^{q+2}}   +  x^{\alpha+1} \otimes  x^{2^{q+2}-2 - \alpha} \otimes z_{2^{q+2}}  )  \\
\\
\hspace{80pt}  
+  x^{\alpha} \otimes  x^{2^{q+2}-2 - \alpha} \otimes Sq^{1}(z_{2^{q+2}} )  \\
\\
=  \alpha ( x^{\alpha} \otimes  x^{2^{q+2}-1 - \alpha} \otimes z_{2^{q+2}}   +  x^{\alpha+1} \otimes  x^{2^{q+2}-2 - \alpha} \otimes z_{2^{q+2}}  )  
\end{array}   \]
since $Sq^{1}(z_{2^{q+2}} )$ is trivial. Therefore
\[   f(x^{2^{q+2}-1} \otimes z_{2^{q+2}})  \in  h({\rm Im}(Sq^{1})  \equiv  0  \]
because $h$ is a $B$-module homomorphism and $Sq^{1}$ is trivial on ${\mathbb Z}/2$.  $\Box$
\end{em}
\end{dummy}
\begin{remark}
\label{2.3}
\begin{em}

When $s = 2,3$ in the situation of Theorem \ref{1.2} there is a map $\alpha :  \Sigma^{\infty} {\mathbb RP}^{\infty} \wedge  {\mathbb RP}^{\infty}  \longrightarrow   \Sigma^{\infty} {\mathbb RP}^{\infty} $ but it is just not equal to $H$! In the loopspace structure of $Q {\mathbb RP}^{\infty}$ form the product minus the two projections to give a map
${\mathbb RP}^{\infty} \times {\mathbb RP}^{\infty} \longrightarrow Q {\mathbb RP}^{\infty}$ which factors through the smash product. The adjoint of this factorisation is $\alpha$. Then the smash product of  two copies of a map of Hopf invariant one $\Sigma^{\infty} S^{2^{s}-1} \longrightarrow  \Sigma^{\infty}{\mathbb RP}^{\infty}$ composed with $\alpha$
is detected by $Sq^{2^{s}}$ on its mapping cone (see \cite{ST82}).

\end{em}
\end{remark}

 \end{document}